\documentclass[a4paper,11pt]{amsart}
\usepackage{geometry}                
\usepackage{graphicx}
\usepackage{amssymb}
\usepackage{xspace}
\usepackage{epstopdf}
\usepackage{tabularx}
\usepackage{enumerate}
\usepackage{floatpag}
\usepackage{geometry}
\usepackage{multirow}
\usepackage{tikz}
\usepackage{longtable}
\usepackage{booktabs}
\usepackage{url}
\usepackage{hyperref}
\usepackage{caption}

\newcommand{\GAP}{{\textsf{GAP}}\xspace}
\newcommand{\Magma}{{\textsc {Magma}}\xspace}

\newcommand{\F}{{\mathbb F}}
\newcommand{\TG}{{\mathsf{TransitiveGroup}}}

\newtheorem{thm}{Theorem}[section]

\theoremstyle{definition}

\DeclareMathOperator{\SL}{\mathrm{GL}}
\DeclareMathOperator{\GL}{\mathrm{GL}}
\DeclareMathOperator{\QD}{\mathrm{QD}}

\DeclareMathOperator{\Alt}{\mathrm{A}}
\DeclareMathOperator{\Sym}{\mathrm{S}}

\newcommand{\Aut}{\mathop{\mathrm{Aut}}}

\newcommand{\trans}{\mathrm{trans}}
\numberwithin{equation}{section}

\renewcommand{\le}{\leqslant}

\newenvironment{mylist}{\begin{list}{}{
\setlength{\parskip}{0mm}
\setlength{\topsep}{2mm}
\setlength{\parsep}{0mm}
\setlength{\itemsep}{0.5mm}
\setlength{\labelwidth}{7mm}
\setlength{\labelsep}{3mm}
\setlength{\itemindent}{0mm}
\setlength{\leftmargin}{12mm}
\setlength{\listparindent}{6mm}
}}{\end{list}}

\title{The transitive groups of degree 48 and some applications}
\author{Derek Holt}
\address{Mathematics Institute, University of Warwick}
\email{D.F.Holt@warwick.ac.uk}

\author{Gordon Royle}
\address{Mathematics and Statistics, University of Western Australia}
\email{gordon.royle@uwa.edu.au}

\author{Gareth Tracey}
\address{Mathematical Institute, University of Oxford}
\email{garethtracey1@gmail.com}

\keywords{transitive group; vertex-transitive graph; census; catalogue;
generator number}
\subjclass[2010]{20B40, 20B35, 05C30}

\begin{document}

\begin{abstract}
The primary purpose of this paper is to report on the successful enumeration
in \Magma of representatives of the $195\,826\,352$ conjugacy classes of
transitive subgroups of the symmetric group $S_{48}$ of degree 48.
In addition, we have determined that 25707 of these groups are minimal
transitive and that 713 of them are elusive.
The minimal transitive examples have been used to enumerate the
vertex-transitive graphs of degree $48$, of which there are
$1\,538\,868\,366$, all but $0.1625\%$ of which arise as Cayley graphs.
We have also found that the largest number of elements required
to generate any of these groups is 10, and we have used this fact to
improve previous general bounds of the third author on the number of elements
required to generate an arbitrary transitive permutation group of a given
degree. The details of the proof of this improved bound will be published
by the third author as a separate paper.
\end{abstract}

\maketitle

\section{Introduction}\label{sec:intro}
Since late in the 19th century, significant effort has been devoted to
compiling catalogues and databases of various types of groups, including
complete lists of (representatives of the conjugacy classes of) the transitive
and primitive subgroups of the symmetric groups of small degree.
For the transitive groups, earlier references include
\cite{GAMillerCollected, Miller96} (with corrections in \cite{MR922391})
for degrees up to $12$, \cite{MR2168238} for degrees up to $31$, 
\cite{MR2455702} for the significantly more difficult case of degree $32$,
and \cite{MR4109707} for degrees $33-47$. (There are $2\,801\,324$ groups of
degree $32$, and a total of $501\,045$ groups of all other degrees up to $47$.)
Apart from the early work of Miller, these lists have been compiled by
computer, using \GAP and \Magma. 

Degree 48 is once again significantly more difficult than earlier degrees,
because there are many more groups and the computations involved need more time
and computer memory. The main purpose of this paper is to report on
the successful enumeration of conjugacy class representatives of the
transitive subgroups of degree $48$, which is the topic of
Section~\ref{sec:main}. There are a total of $195\,826\,352$ of these subgroups.
The computations were carried out in \Magma and required a total
of about one year of cpu-time.
The complete list of these subgroups is available in \Magma using an
optional database that can be downloaded by users from the \Magma website.
Although we carried out these computations serially on a single processor and,
due to various logistical problems, they took more than two years of real time
to complete, they are intrisically extremely parallelisable: about $98\%$ of the
cpu-time was for imprimitive groups with blocks of size $2$ and, as we shall
explain shortly, that case splits into $25000$ independent calculations.

We anticipate that it would be feasible to extend the catalogues up to degree
$63$, but that degree $64$ will remain out of range for the foreseeable future.

\begin{table}

\begin{tabular}{|r|rrr|}
\hline
MinBlockSize & TransGps & MinTrans & Elusive\\
\hline
2 &      192327620 & 15046 & 108\\
3 &        3397563 & 10625 & 590\\
4 &          94121 &    36 &   3\\
6 &           5011 &     0 &  12\\
8 &           1275 &     0 &   0\\
12 &           103 &     0 &   0\\
16 &           646 &     0 &   0\\
24 &             9 &     0 &   0\\
Primitive &      4 &     0 &   0\\
\hline
Total &  195826352 & 25707 & 713\\
\hline
\end{tabular}
\vspace{5pt}

\caption{Numbers of transitive, minimal transitive, and elusive
 groups of degree 48}
\label{tab:summary}
\end{table}

We subsequently used our catalogue to identify those groups of degree $48$ that
are minimal transitive (that is, they have no proper transitive subgroups)
and those that are \emph{elusive} (that is, they contain no fixed-point-free
elements of prime order).  We shall report on this in Section~\ref{sec:minel}.
The various counts of groups involved are summarised in Table~\ref{tab:summary},
where the imprimitive groups have been counted according to the smallest
size of a block of imprimitivity.

In Section~\ref{sec:graphs}, we describe the computation of the vertex-transitive
graphs of order 48, along with some associated data.

We denote the smallest size of a generating set of a group $G$ by $d(G)$.
We have established by routine computations, that $d(G) \leqslant 10$ for all
transitive groups of degree $48$.  In fact, the only examples with $d(G)=10$
have minimal block size $3$ and have $10$-generator transitive groups of degree
$32$ as quotients. The groups with block size $2$ all satisfy $d(G) \leqslant 9$.

As we shall explain in Section~\ref{sec:gennos}, these bounds on $d(G)$
have enabled the third author to remove the exceptional cases of the general
bound on $d(G)$ for transitive permutation group of degree $n$ that he
established in~\cite{MR3812195}, and thereby to complete the proof of the
following result, where logarithms are to the base $2$.
\begin{thm}\label{thm:transgens}
Let $G$ be a transitive permutation group of degree $n$. Then
$$d(G)\leqslant \left\lfloor\frac{cn}{\sqrt{\log{n}}}\right\rfloor$$
where $c:=\frac{\sqrt{3}}{2}$.
\end{thm}

(It was proved by Lucchini in \cite{Luc} that this result holds for some
unspecified constant $c$.)
Since the proof of this result involves some lengthy case-by-case analyses,
we shall just summarise it in Section~\ref{sec:gennos} of this paper, and
the details will be published separately by the third author~\cite{Tracey21}.

In a related application, the third author is now able to improve a
previously unpublished result bounding the constant $d$ in the result
proved in \cite{LMM} that $d(G) \le d\log  n/\sqrt{\log \log  n}$ 
for primitive subgroups $G$ of $S_n$.

\vspace{0.2cm}
\noindent\textbf{Notation:} For a finite group $G$, we will write $\Phi(G)$, $R(G)$, $[G,G]$ for the Frattini subgroup, soluble radical, and derived subgroup of $G$, respectively. We will mostly use the notation from \cite{WebAtlas} for group names, although we simply write $n$ for the cyclic group of order $n$ when there is no danger of confusion.

\section{Computing the transitive groups of degree $48$}
\label{sec:main}
The primitive permutation groups are known up to degree $4095$\footnote{
This has recently been extended to degree $8191$ by Ben Stratford, a
student of the first author.}
(see \cite{MR2845584}) and are incorporated into the databases of both
\Magma and \GAP, and so we need only consider the imprimitive groups. 
By definition, if a group $G$ acting transitively on the set $\Omega$ of size
$n$ is not primitive, then there is at least one partition of $\Omega$ into a
block system $\mathcal{B}$ such that $G$ permutes the blocks of $\mathcal{B}$.
If we let $G^\mathcal{B}$ denote the action of $G$ on the blocks
(``the top group'') and if $\mathcal{B}$ has $n/k$ blocks of size $k$, 
then $G^\mathcal{B}$ is a transitive permutation group of degree $n/k$. 
We say that the block system $\mathcal{B}$ is \emph{minimal} if $k$ is
minimal among block systems with $k>1$.
Then we can associate to each group $G$ a set of pairs of the form
\[
\{(k, G^\mathcal{B}) : \mathcal{B}\text{ is a minimal block system for } G \text{ with blocks of size } k\}.
\]

If this set contains more than one pair (imprimitive groups may of course have
more than one minimal block system), then we wish to distinguish just one of
them. Thus we define the \emph{signature} of an imprimitive permutation group
to be the \emph{lexicographically least} pair $(k,G^\mathcal{B})$ associated
with $G$, where the second component is indexed according to its order in the
list of transitive groups of degree $n/k$ already in {\sc Magma}. 
But note that it can happen that two different minimal block systems of $G$
define the same signature.

We separate the computation into parts, with each part constructing only the
groups with a particular signature. Given an integer $k$ such that $1 < k < n$,
and a transitive group $H$ of degree $n/k$, the wreath product
$S_k \wr H$ contains (a conjugate of) every transitive group of degree $n$
with signature $(k,H)$. So these groups can all be found by exploring the
subgroup lattice of $S_k \wr H$ (although there are complications arising from
the fact that we want representatives of subgroups up to conjugacy in $S_n$.)

A naive approach to the problem for a fixed $k$, is to deal with all
candidates $H$ simultaneously,  by starting with $S_k \wr S_{n/k}$,
and repeatedly using the \texttt{MaximalSubgroups} command of \Magma,
thereby traversing the subgroup
lattice downwards and in a breadth-first fashion, pruning each branch of the
search as soon as it produces groups with signature differing from $H$,
while using conjugacy tests to avoid duplication.
(We also have to eliminate duplicates arising from a group preserving
more than one minimal block system with blocks of size $k$.)
This was successfully applied in all cases to the transitive groups of
degrees $33-47$, and we refer the reader to \cite[Section 2]{MR4109707}
for further details. In degree $48$, we successfully applied this
method to groups with signatures $(k,H)$ with $k \geqslant 6$; that is for $k=6,8,12,
16$ and $24$. The cpu-times in these cases were of order 10 hours,
30 minutes, 3 minutes, 70 minutes, and a few seconds, respectively.

For the examples with $k=2,3$ and $4$, we used the methods described in
\cite[Section 3]{MR4109707} (and also in \cite[Section 2.2]{MR2455702} for
$k=2$).

The methods for $k=3$ were essentially the same as in degree $36$, but
considerably more time-consuming.
We have $G \le S_3 \wr S_{16} \cong C_3 \wr (C_2 \wr S_{16})$.
Let $\rho$ be the induced projection of $G$ onto $C_2 \wr S_{16}$. Then,
since we are assuming that $G$ is transitive, $\rho(G)$ must project onto a
transitive subgroup of $S_{16}$, and the existing catalogues contain the
$1954$ possibilities for this projection. Furthermore, either
\begin{enumerate}
\item[(i)] $\rho(G)$ is a transitive subgroup of $S_{32}$, in which case we
can use the existing catalogue
$2\,801\,324$ as a list of candidates for $\rho(G)$; or
\item[(ii)] $\rho(G)$ is an intransitive groups of degree $32$ that projects
onto a transitive subgroup of $S_{16}$. In that case, it is not hard to show
that $G$ must be conjugate to the natural complement of the base group of
$C_2 \wr H$, where $H$ is one of the 1954 transitive groups of degree $16$.
(We also checked this computationally.)
\end{enumerate}

We enumerated the groups with $k=3$ by considering each of the
$2\,801\,324 + 1954$ possibilities for $\rho(G)$ in turn.
This involves a cohomology computation, which is analogous to that for
the case $k=2$, which we shall discuss below.
The computation for $k=3$ took a total of about 104 hours of cpu-time.
Of the $3\,397\,563$ groups on this list, $\rho(G)$ is transitive for all
except $55\,715$. 
Similarly, for $k=4$, we used the same techniques as in degrees $36$ and $40$,
and the total cpu-time was about 9 hours.

The vast majority of the computational work was for the case $k=2$, and
we shall briefly recall how we proceed in this case.
We have $G \le W := C_2 \wr H$, where $H := G^\mathcal{B}$ is one of the groups
in the known list of $25\,000$ transitive groups of degree $24$.
Again we calculate those groups with signature $(2,H)$ for each
individual group $H$, and the $25\,000$ calculations involved are
independent and could in principal be done in parallel.

Let $K \cong  C_2^{24}$ be the kernel of the action of $W$ on $\mathcal{B}$.
Then we can regard $K$ as a module for $H$ over the field
$\F_2$ of order $2$, and $M := G \cap K$ is an
$\F_2H$-submodule. We can use the \Magma commands \texttt{GModule}
and \texttt{Submodules} to find all such submodules.
In fact, since we are looking for representatives of the conjugacy
classes of transitive subgroups of $W$, we only want one representative 
of the conjugation action of $N := C_2 \wr N_{S_{24}}(H)$
on the set all $\F_2H$-submodules $M$ of $K$, and we
use the \Magma command \texttt{IsConjugate} to find such representatives.

Now, for each such pair $(H,M)$, the transitive groups $G$ with
$H = G^\mathcal{B}$ and $M = G \cap K$ correspond to complements
of $K/M$ in $H/M$, and the $H$-conjugacy classes of such complements
correspond to elements of the cohomology group $H^1(H,K/M)$, which
can be computed in \Magma.

We also need to test these groups $G$ for conjugacy under the action of
$N_N(M)$. In some cases when $H^1(H,K/M)$ is reasonably small, this can be done
in straightforward fashion using {\sc Magma}'s \texttt{IsConjugate} function.
But in many cases this was not feasible, and we had to use the method
using an induced action on the cohomology group that is described in detail
in \cite[Section 2.2]{MR2455702}. Finally, for each $G$ that we find, we need
to find all block systems with block size $2$ preserved by $G$, so that we can
eliminate occurrences of groups that are conjugate in $S_n$ but arise either for
distinct pairs $(H,M)$ or more than once for the same pair. Again
we refer the reader to  \cite[Section 2.2]{MR2455702} for further details.

Here are some statistical details concerning some of these calculations.

\begin{itemize}
\item The numbers of groups arising from the $25000$ candidates for the top
group $H$ ranges from $3$ to $3\,642\,186$, with average $7693$ and median
$778$.  This number is less than $10000$ for  more than $90\%$ of the
top groups.  For the majority of these groups $H$, the computations were fast.
For examples, for the groups $H = \TG(24,k)$ with $20000 < k \le 25000$,
($20\%$ of the groups $H$) the total cpu-time was about $92.4$ hours
(just over $1\%$ of the total) and the total number of groups $G$ that arise is
$2\,963\,853$ (about $1.5\%$ of the total).
\item
The highest dimension of a cohomology group $H^1(H,M)$ was $26$. In that case,
$|H^1(H,M)| = 2^{26} = 67\,108\,864$, where elements of $H^1(H,M)$ are
represented by 26 binary digits. This is important because of
the orbit computation on the elements of $H^1(H,M)$. If an example of much
higher dimension than this had been encountered (which we might expect to
be the case for a corresponding attempt to find the transitive groups of
degree 64), then this orbit computation might not have been feasible.
This occurred with $H = \TG(24,4010)$ and $|M|=2^{12}$,
and the pair $(H,M)$ gave rise to 201\,792 groups $G$.

\item
The case $H = \TG(24,11363)$ resulted in the largest number of groups $G$,
namely $3\,642\,186$. There were 240 possibilities for $M$. This case took
about 34 hours of cpu-time, using about 73GBytes RAM.

\item
The pair $(H,M)$ that resulted in the most groups $G$, namely $1\,054\,720$,
arose with $H = \TG(24,13329)$ and $|M| = 2^7$. Although $H^1(H,M)$ had
dimension only 22 in this case, there were many more orbits of the action
than in the case with dimension 26 discussed above.
\end{itemize}

\section{Minimal transitive and elusive groups}
\label{sec:minel}
\subsection{Minimal transitive groups}
For many applications that involve considering all possible transitive actions
of a certain degree, it is sufficient to consider only the
\emph{minimal transitive groups} i.e., transitive groups with no proper
transitive subgroups. (One example of this was discussed in
\cite[Section 5]{MR4109707}, where all vertex-transitive graphs of
degrees 33--47 are constructed.) Testing if a transitive group is minimal can
be done by finding all of its maximal subgroups and verifying that none are
transitive. As most of the groups are not minimal transitive, it proves useful
in practice to first construct some random subgroups in an attempt to find a
transitive proper subgroup, only undertaking the more expensive step of finding
all maximal subgroups if this fails.

There are a total of 25707 minimal transitive groups, all of which have minimal
blocks of sizes $2$, $3$ or $4$ --- the exact numbers of minimal transitive
groups with smallest blocks of each size are given in Table~\ref{tab:summary}.

\subsection{Elusive groups}

One of the major reasons to construct catalogues of combinatorial objects is
to gather evidence relating to conjectures or other open questions.
Even if a newly-constructed catalogue does not directly contain a
counterexample to a conjecture (thereby immediately resolving it), it can be
useful in refining a researcher's intuition regarding both the typical and
extremal objects in the catalogue. 

A permutation group $G$ is called {\em elusive} if it contains no
\emph{fixed-point-free} elements (i.e., \emph{derangements}) of prime order.
Elusive groups are interesting because of their connection to
Maru\v{s}i\v{c}'s \emph{Polycirculant Conjecture} \cite{PCC} which asserts that the
automorphism group of a vertex-transitive digraph is \emph{never elusive}.
In principle, a positive resolution of the polycirculant conjecture may
simplify the construction and analysis of vertex-transitive graphs and
digraphs, as it would then be possible to assume the presence of an
automorphism with $n/p$ cycles of length $p$ for some prime $p$. 
Early catalogues of vertex-transitive graphs often used ad hoc arguments to
show that all transitive groups of the specific degrees under consideration
have a suitable derangement of prime order.

A permutation group $G$ is called \emph{$2$-closed} if there is no group
properly containing $G$ with the same orbitals as $G$. The automorphism group
of a vertex-transitive digraph is necessarily $2$-closed, because it is
already the maximal group (by inclusion) that fixes the set of arcs of the
digraph, which is a union of some of the orbitals. The conjecture can thus be
strengthened to the assertion that there are no elusive $2$-closed transitive
groups, as proposed by Klin and Maru\v{s}i\v{c} at the 15th British Combinatorial
Conference \cite{Klin}.

One might hope  that there are simply no elusive groups at all, in which case
both conjectures would hold vacuously, but in fact there are a number of
sporadic examples of elusive groups and a handful of infinite families. However all 
the known elusive groups are not $2$-closed, so do not provide counterexamples
for either conjecture.  

It is relatively easy to test the groups for the property of being elusive by
checking to see if any of the conjugacy class representatives are derangements
of prime order. For the larger groups, it is often faster to first generate
some number of randomly selected elements inside each of the Sylow subgroups
in the hope of stumbling on a suitable derangement without the cost of
computing all the conjugacy classes.

The results of this computation reveal that there are 713 elusive groups of degree 48, with orders 
ranging from $5184$ to $806\,215\,680\,000$.
The numbers of elusive groups of degree 48 with each minimal
blocksize are given in Table~\ref{tab:summary}. If an elusive 
group has minimal blocks of different sizes (say 2 and 3), then it is grouped
and counted according the smaller of the sizes.


Of these groups $700$ have a unique minimal normal subgroup, and while each of the remaining $13$ groups has multiple minimal normal subgroups, these minimal normal subgroups are conjugate in $S_{48}$. Therefore we can partition the elusive groups according to the unique conjugacy class of their minimal normal subgroup(s).

Collectively, the $713$ elusive groups share just $7$ pairwise non-conjugate minimal normal subgroups. Table~\ref{fig:mns} shows the different minimal normal subgroups that occur and the number of elusive groups with that particular minimal normal subgroup. In addition, it gives the order of the normalizer (in $S_{48}$) of that subgroup, while the final column shows the possible minimal block sizes that occur for that minimal normal subgroup. All but one of the possible minimal normal subgroups are elementary abelian, but two non-conjugate (but obviously isomorphic) groups of orders $2^8$ and $3^8$ occur. For example, the first row shows that an elusive group with minimal normal subgroup $C_3^{4}$ either has minimal blocks of size $2$ (only) or minimal blocks of sizes both $2$ and $3$.

\begin{table}
\begin{center}
\renewcommand{\arraystretch}{1.1}
\begin{tabular}{lcccc}
\toprule
Group & $\vert$Normalizer$\vert$ & Frequency & Min. Blocks\\
\midrule
$C_3^4$ & $2^{17} \cdot 3^{18}$ & 22 & $\{2\}$ or $\{2,3\}$ \\
$C_2^8$ & $2^{28} \cdot 3^3$ &3 & $\{4\}$ \\
$C_2^8$ & $2^{43} \cdot 3^5$ & 75 & $\{2\}$ \\
$C_3^8$ & $2^{19} \cdot 3^{21}$ & 575 & $\{3\}$\\
$C_3^8$ & $2^{11} \cdot 3^{17} \cdot 7$ & 2 & $\{3\}$\\
$C_2^{16}$ & $2^{39} \cdot 3^{10} \cdot 5 \cdot 7$&24& $\{2\}$\\
$\Alt(6)^4$ & $2^{23} \cdot 3^{9} \cdot 5^4$ & $12$ & $\{6\}$\\
\bottomrule
\end{tabular}
\caption{Minimal normal subgroups of elusive groups of degree 48}
\label{fig:mns}
\end{center}
\end{table} 

With this many elusive groups, and no obvious way to get a
compact description, it would seem unlikely that the Polycirculant
Conjecture can be proved by first classifying elusive groups.

\section{Vertex-transitive graphs of order 48}\label{sec:graphs}

The class of vertex-transitive graphs plays a central role in algebraic graph theory, often 
providing extremal cases or illuminating examples in the study of many graphical properties.
Although it is not strictly necessary to have a complete list of the transitive groups of degree $d$ in
order to compute a complete list of vertex-transitive graphs of order $d$, it is conceptually simple
to compute all the vertex-transitive graphs from the transitive groups.  

For notational convenience, we say that a graph $\Gamma$ is $G$-vertex-transitive (or just
$G$-transitive) if $G \leqslant \Aut(\Gamma)$ and $G$ acts transitively on $V(\Gamma)$.
Given a list of all the transitive groups of some fixed degree, in principle it suffices
to consider each group $G$ in turn, construct all the $G$-transitive
groups, and then merge the lists from the different groups, removing all but one
isomorphic copy of each graph.

As stated, this naive algorithm would do far too much work, constructing large numbers of isomorphic 
copies of most of the graphs. However we can reduce this work in two ways. First we can restrict
our attention to the \emph{minimal transitive} groups, because if $H \leqslant G$ and $H$ is transitive,
then any $G$-transitive group is $H$-transitive. Secondly, we can do some work to avoid 
constructing graphs that are \emph{obviously isomorphic} to ones that 
have been, or will be, constructed elsewhere.

The transitive groups of degree $48$ and order $48$ are necessarily minimal transitive, and we deal with these separately from
the larger minimal transitive groups.  This separates out the vertex-transitive graphs with a regular subgroup of automorphisms, i.e., the
Cayley graphs, from the remainder. This is common practice, because there are a number of interesting questions
and conjectures where the distinction between Cayley graphs and non-Cayley graphs appears to be subtle but significant.

\subsection{Cayley graphs}

Given a group $G$, and a set of group elements  $C \subseteq G$ such that ${\bf 1}_G \notin C$ and $C^{-1} = C$, 
the \emph{Cayley graph} ${\rm Cay}(G,C)$ is the graph defined as follows:
\begin{align*}
V({\rm Cay}(G,C)) &= G,\\
E({\rm Cay}(G,C)) &= \{(g,cg) \mid g \in G, c \in C\}.
\end{align*}

It is immediate that $G$ is a regular subgroup (acting by right-multiplication) of the 
automorphism group of  ${\rm Cay}(G,C)$, and it is well-known that any vertex-transitive graph 
whose automorphism group has a regular subgroup $R$ is a Cayley graph for $R$. 
The set $C$ is called the \emph{connection set} for the Cayley graph and it is precisely 
the neighbourhood of the vertex ${\bf 1}_G$. (The conditions on $C$ are simply to ensure
that the resulting graphs are undirected and loopless.)

While all Cayley graphs are vertex transitive, not all vertex-transitive graphs are Cayley graphs, with the
canonical example here being the Petersen graph.  For small orders, the vast majority of 
vertex-transitive graphs are Cayley graphs, but it is not known if this holds in general. In other words,
is it true that the proportion of vertex-transitive graphs of order at most $n$ that are Cayley graphs tends
to $1$ as $n$ increases? 

If we define
\(
\Omega = \{ \{g,g^{-1}\} \mid g \in G, g \ne {\bf 1}_G \}
\)
then every subset of $\Omega$ determines the connection set for some Cayley graph of $G$ and vice versa. If  $G$ has $a$ involutions, and $b$ non-identity element-inverse pairs, then $|\Omega| = a+b$, and so there are exactly $2^{a+b}$ Cayley graphs for $G$.
As we are almost always only interested in isomorphism classes of graphs, we need to remove, or preferably never construct, all but one representative of each isomorphism class. 
There is one obvious source of isomorphisms, namely those  arising from the automorphism group of $G$. More precisely, if $\sigma \in \Aut(G)$ then ${\rm Cay}(G,C)$ is isomorphic to ${\rm Cay}(G,C^\sigma)$.  So if we define a \emph{Cayley set} to be an orbit of $\Aut(G)$ acting on the set of subsets of $\Omega$, then it suffices to consider just one connection set from each Cayley set.  

In practice, we fix an arbitrary order on $\Omega$, use {\sffamily GAP} to compute the action of $\Aut(G)$ on $\Omega$, and then a simple orderly-style algorithm to compute the lexicographically least representative of each Cayley set. This makes heavy use of Steve Linton's \verb+SmallestImageSet+ \cite{Linton} to ensure that only lexicographically-least subsets of $\Omega$ are considered at every stage.  The theoretical number of Cayley sets for $G$ can be determined by calculating the cycle index polynomial of $\Aut(G)$ acting on $\Omega$ (using the undocumented  \verb+CycleIndexPolynomial+ function in {\sc Magma}) and applying P\'{o}lya's Enumeration Theorem. For each of the 52 groups examined, the actual number of Cayley sets constructed by the orderly algorithm matches the theoretical number, giving us a high degree of confidence in this stage of the computation. 

Using the list of (representatives of) Cayley sets we next construct the corresponding list of Cayley graphs. Although this list is free of isomorphisms induced by the action of $\Aut(G)$ on the connection sets, there can be additional isomorphisms. Therefore we filter the list of Cayley graphs for each group, removing any graph that is isomorphic to an earlier graph in the list. Sometimes there are no isomorphisms except the ones induced by  $\Aut(G)$, and so the final filtering step does not remove any graphs. Groups with this property are called {\sffamily CI}-groups and there is a substantial literature on the still-open question of characterizing {\sffamily CI}-groups.

There are 52 groups of order $48$ and the results of the Cayley graph computations for those groups are given in Table~\ref{fig:cayley48}. The groups are numbered from $1$ to $52$ according to their order in the small group libraries of  \Magma and {\sffamily GAP} (the groups are in the same order in each library). The  group structure is the description returned by the {\sffamily GAP} command \verb+StructureDescription+, the values $a$ and $b$ are the number of involutions and the number of non-identity element-inverse pairs respectively. The column labelled $|\mathrm{Aut}|$ lists the order of the automorphism group of the group. In most cases, the value $2^{a+b}/|\mathrm{Aut}|$ is an approximation (an underestimate) for the number of Cayley sets for that group. Where an entry in the  $|\mathrm{Aut}|$ column is  marked with an asterisk (such as $^*192$ for group number $8$), this indicates a group where $\Aut(G)$ does not act faithfully on $\Omega$. These groups are characterized by the existence of a group automorphism $\sigma \in \Aut(G)$ such that $g^\sigma \in \{g, g^{-1}\}$ for each element $g$. It is known that such a group is either an \emph{abelian group} (where the inverse map is a group automorphism) or a \emph{generalised dicyclic group} (Watkins \cite{Watkins}).
In each of the cases indicated in Table~\ref{fig:cayley48}, the kernel of the action of  $\Aut(G)$ on $\Omega$ has order $2$, and so these groups yield approximately $2^{a+b+1}/|\mathrm{Aut}|$ Cayley sets.

The column ``Cayley Sets'' gives the exact number of Cayley sets obtained from Polya's Enumeration Theorem, while the final column ``Cayley Graphs'' gives the actual number of pairwise non-isomorphic Cayley graphs. This last step is computationally non-trivial because of the sheer size of some of the lists. For example, there are more than 360 million Cayley graphs for the most prolific group ($C_2 \times D_{24}$).

A group is a CI-group if and only if the number of Cayley graphs is equal to the number of Cayley sets, so the table shows that the abelian group $C_{6} \times C_{2} \times C_{2} \times C_{2}$ is the only CI-group of order 48.

{\footnotesize
\begin{center}
\def\arraystretch{0.84}
\begin{longtable}[p]{ccrrrrr}
\caption{Cayley graphs on 48 vertices}\label{fig:cayley48}\\
\toprule
No. & Structure & $a$ & $b$ & $|\mathrm{Aut}|$ & Cayley Sets & Cayley Graphs\\
\midrule
1&$C_{3} : C_{16}$&1&23&48&496512&489376\\
2&$C_{48}$&1&23&$^*$16&2151936&2122944\\
3&$(C_{4} \times C_{4}) : C_{3}$&3&22&384&104224&103726\\
4&$C_{8} \times S_3$&7&20&48&3752448&3516448\\
5&$C_{24} : C_{2}$&7&20&48&3561216&3337160\\
6&$C_{24} : C_{2}$&13&17&96&13641984&11880240\\
7&$D_{48}$&25&11&192&364086016&360716112\\
8&$C_{3} : Q_{16}$&1&23&$^*$192&275712&255696\\
9&$C_{2} \times (C_{3} : C_{8})$&3&22&96&647168&597648\\
10&$(C_{3} : C_{8}) : C_{2}$&3&22&96&586752&553168\\
\midrule
11&$C_{4} \times (C_{3} : C_{4})$&3&22&192&454176&370704\\
12&$(C_{3} : C_{4}) : C_{4}$&3&22&96&893952&611760\\
13&$C_{12} : C_{4}$&3&22&$^*$192&586752&484944\\
14&$(C_{12} \times C_{2}) : C_{2}$&15&16&96&28924416&23139848\\
15&$(C_{3} \times D_8) : C_{2}$&17&15&96&47661696&45855520\\
16&$(C_{3} : Q_8) : C_{2}$&5&21&96&1102464&967024\\
17&$(C_{3} \times Q_8) : C_{2}$&13&17&96&12473472&11952272\\
18&$C_{3} : Q_{16}$&1&23&96&358272&308400\\
19&$(C_{6} \times C_{2}) : C_{4}$&7&20&192&1429024&1049296\\
20&$C_{12} \times C_{4}$&3&22&$^*$192&452032&431808\\
\midrule
21&$C_{3} \times ((C_{4} \times C_{2}) : C_{2})$&7&20&64&3373440&2876192\\
22&$C_{3} \times (C_{4} : C_{4})$&3&22&64&1081344&904064\\
23&$C_{24} \times C_{2}$&3&22&$^*$32&2336768&2183232\\
24&$C_{3} \times (C_{8} : C_{2})$&3&22&32&1603584&1499008\\
25&$C_{3} \times D_{16}$&9&19&64&5257920&5124112\\
26&$C_{3} \times \QD_{16}$&5&21&32&2776320&2592224\\
27&$C_{3} \times Q_{16}$&1&23&64&447168&423232\\
28&$C_{2} . S_4 = \SL(2,3) . C_{2}$&1&23&48&436864&431120\\
29&$\GL(2,3)$&13&17&48&22758528&22589392\\
30&$A_4 : C_{4}$&7&20&48&3343616&3230252\\
\midrule
31&$C_{4} \times A_4$&7&20&48&3203072&3142848\\
32&$C_{2} \times \SL(2,3)$&3&22&48&808832&791338\\
33&$((C_{4} \times C_{2}) : C_{2}) : C_{3}$&7&20&48&2992000&2967742\\
34&$C_{2} \times (C_{3} : Q_8)$&3&22&$^*$384&386784&279616\\
35&$C_{2} \times C_{4} \times S_3$&15&16&192&16697472&14159528\\
36&$C_{2} \times D_{24}$&27&10&384&373069248&362458536\\
37&$(C_{12} \times C_{2}) : C_{2}$&15&16&96&29701632&22395608\\
38&$D_8 \times S_3$&23&12&96&383280384&349815008\\
39&$(C_{4} \times S_3) : C_{2}$&11&18&96&9219840&5452760\\
40&$Q_8 \times S_3$&7&20&288&1108224&823752\\
\midrule
41&$(C_{4} \times S_3) : C_{2}$&19&14&288&34347520&30759480\\
42&$C_{2} \times C_{2} \times (C_{3} : C_{4})$&7&20&$^*$1152&428256&357264\\
43&$C_{2} \times ((C_{6} \times C_{2}) : C_{2})$&19&14&192&54719616&47641176\\
44&$C_{12} \times C_{2} \times C_{2}$&7&20&$^*$384&967808&806656\\
45&$C_{6} \times D_8$&11&18&128&6653184&5654768\\
46&$C_{6} \times Q_8$&3&22&384&246912&207744\\
47&$C_{3} \times ((C_{4} \times C_{2}) : C_{2})$&7&20&96&2397184&1865216\\
48&$C_{2} \times S_4$&19&14&48&182656512&177923704\\
49&$C_{2} \times C_{2} \times A_4$&15&16&144&15715936&15306700\\
50&$(C_{2} \times C_{2} \times C_{2} \times C_{2}) : C_{3}$&15&16&5760&413344&411248\\
51&$C_{2} \times C_{2} \times C_{2} \times S_3$&31&8&8064 & 72984704&71039696\\
52&$C_{6} \times C_{2} \times C_{2} \times C_{2}$&15&16&$^*$40320& 160208 & 160208\\
\bottomrule
\end{longtable}

\end{center}}

The total number of Cayley graphs, after removing isomorphs both \emph{within} and \emph{between} the 52 lists of graphs is $1\,536\,366\,616$, approximately $1.54$ billion Cayley graphs of order $48$.

\subsection{Non-Cayley graphs}

The automorphism group of a vertex-transitive graph that is not a Cayley graph 
is a transitive group that has no regular subgroups. Thus the simplest approach 
is to consider each minimal transitive group $G$ in turn, compute
all the $G$-transitive graphs and then remove both unwanted isomorphs
\emph{and} any Cayley graphs that have accidentally been constructed along
the way.


If $G$ is a transitive group then its \emph{orbitals} are defined to be the orbits of $G$
on $V(G) \times V(G)$.  If $\mathcal{O} = (x,y)^G$ is an orbital of $G$, then $(y,x)^G$ is 
also an orbital of $G$, called the \emph{paired orbital} of $\mathcal{O}$. A graph is $G$-transitive
if and only if its edge-set is the union of pairs of orbitals of $G$ (identifying an edge $xy$ with a 
pair of oppositely-directed arcs $\{(x,y), (y,x)\}$).

An \emph{orbital graph} of $G$ is a graph whose edge-set is the union of a single pair of
orbitals. Let $G'$ be the intersection of the automorphism groups of all the orbital
graphs of $G$. (This is essentially an undirected version of the 2-closure of a group, sometimes
called the \emph{strong $2$-closure} of the group.)
Then any $G$-transitive graph is $G'$-transitive and so it is only
necessary to process $G'$. In itself, this does not reduce the amount of work 
required because, by definition, $G$ and $G'$ have precisely the same 
pair-closed sets of orbitals to consider.  However it is often the case
that $G'$ and $H'$ are conjugate even when $G$ and $H$ are not. Therefore
by constructing the strong 2-closure of all of the minimal transitive groups
and then throwing out all-but-one from each conjugacy class, we end up with a much 
smaller list of groups to process. This list can be even further reduced by 
noting that $G'$ sometimes has a regular subgroup although $G$, again by definition, does not.
In this situation, every $G'$-transitive graph is a Cayley graph, 
and as these have already been constructed, there is no need to process $G'$.

The minimal transitive groups of order greater than 48 collectively 
have 840 pairwise non-conjugate strong 2-closures. It is easy to verify that
a small group has no regular subgroups before processing it further, but for
the larger groups this becomes too time-consuming. However the larger
groups tend to have few orbitals, and so it is easy to construct
all of the transitive graphs for these groups. The resulting list of graphs
then contains all of the non-Cayley graphs, but also many Cayley graphs that
must be removed. Due to the sheer size of the computation, this is a 
rather lengthy and somewhat intricate process, but on completion we end up
with $2\,501\,750$ non-Cayley graphs, of which $2\,501\,630$ are connected. 
Hence the total number of vertex-transitive
graphs on 48 vertices is $1\,538\,868\,366$, of which 
just $0.1625\%$ are not Cayley graphs.

\subsection{Edge-transitive and half-arc transitive graphs}

A graph is called \emph{edge-transitive} if its automorphism group is transitive on edges (i.e., unordered pairs of adjacent vertices) and \emph{arc-transitive} if it is transitive on its arcs (i.e., ordered pairs of adjacent vertices). An edge-transitive graph might also be vertex-transitive, but there are edge-transitive graphs that are not vertex transitive, in fact some that are not even regular. Conder \& Verret \cite{MR4049846} have computed all of the edge-transitive graphs on up to 47 vertices, separately finding those that are vertex transitive, and those that are not. 

As a result of the computations reported above, we can go one step further and find the edge-transitive graphs of order $48$ that are also vertex-transitive. Thus from the $1.54$ billion vertex-transitive graphs of order $48$, we extracted $189$ edge-transitive graphs, of which $115$ are connected, $115$ (sic) are twin-free (\emph{twins} are vertices with the same neighbourhood) and just $71$ are both connected and twin-free. 

We can also extract a few more interesting graphs from our lists.  A graph is called \emph{half-arc transitive} (or just \emph{half-transitive}) if it is vertex transitive and edge transitive, but not arc transitive.
The most famous, and smallest, such
graph is the $4$-regular graph on $27$ vertices known as the \emph{Doyle-Holt graph} after its independent discoverers Doyle \cite{Doyle} (originally in an unpublished Masters Thesis at Harvard in 1976) and Holt \cite{MR615008} in 1981.

The data tabulated in Conder \& Verret indicate that there is a single half-arc-transitive graph on $27$ vertices (degree $4$), $2$ on $36$ vertices (of degrees $8$ and $12$), $2$ on $39$ vertices (degrees $4$ and $8$) and $3$ on $40$ vertices (all of degree $8$). To this we can add another $4$ half-arc-transitive graphs on $48$ vertices (all of degree $8$).

All four of these $8$-regular half-arc-transitive graphs are Cayley graphs for at least one group of order $48$, with the groups occurring being 
$(C_{4} \times C_{4}) : C_{3}$ (Group 3 from the list above), $A_4 : C_{4}$ (group 30), $C_{2} \times C_{2} \times A_4$ (group 49) and $(C_{2} \times C_{2} \times C_{2} \times C_{2}) : C_{3}$ (group 50). 

\section{Maximal generating number of transitive groups of degree 48}
\label{sec:gennos}
For an arbitrary group $G$, let $d(G)$ be the minimal size of a generating set
of $G$. As we saw earlier, for most of the transitive groups $G$ of
degree $48$ that are imprimitive with block size $3$, the quotient group
$\bar{G} :=G/O_3(G)$ of $G$ is naturally isomorphic to a transitive group of
degree $32$, There are five such groups with $d(\bar{G}) = 10$, namely
$\mathsf{TransitiveGroup}(32,i)$ for $i \in \{
    1422821, 1422822, 1514676, 2424619, 2224558 \}$, and it turned out that
there are also five corresponding groups $G$, also with $d(G) = 10$.
A lengthy but routine computation showed that these are the only transitive
groups of degree 48 with $d(G) > 9$.

Among the groups $G$ with minimal block size $2$, there are 11 groups with
$d(G) = 9$, and these have signatures $(2,H)$, where
$H = \mathsf{TransitiveGroup}(24,i)$ with $i \in \{9169, 21182, 23560\}$.
The maximum value of $d(G)$ among primitive groups and groups with minimal
block size at least $4$ is $6$, which arises with block sizes 4 and 6 only.

In \cite{MR3812195}, the problem of finding numerical upper bounds for $d(G)$ for an arbitrary transitive permutation group $G$ of degree $n$ is considered. It had already been proved in \cite{Luc} that $d(G)$ is at most $\frac{cn}{\sqrt{\log n}}$ in this case, where $c$ is an unspecified absolute constant. This bound is shown to be asymptotically best possible in \cite{KN} (that is, there exists constants $c_1$, $c_2$, and an infinite family $(G_{n_i})_{i=1}^{\infty}$ of transitive groups of degree $n_i$, with $c_1\le \frac{n_i}{d(G_{n_i})\sqrt{\log n_i}}\le c_2$ for all $i$).

In \cite{MR3812195} it is proved that, apart from a finite list of possible exceptions, the bound $d(G)\le \left\lfloor \frac{cn}{\sqrt{\log n}}\right\rfloor$ holds, where $c:=\frac{\sqrt{3}}{2}$ (and logarithms are to the base $2$).
This bound is best possible in the sense that
$d(G)= \frac{\sqrt{3}n}{2\sqrt{\log n}}=4$
when $G=D_8\circ D_8<\Sym_8$ and $n=8$, although  it seems likely that there
are better bounds that hold for sufficiently large $n$.

The information in the first paragraph above concerning generator numbers in
transitive groups of degree $48$ has helped the third author to complete the
proof of Theorem \ref{thm:transgens} thereby dispensing with the finite list of
exceptions. There are, however, a number of other steps in this proof, some of
which involve lengthy case-by-case analyses. For this reason, we will just
outline the general strategy of the proof in this paper, and the details
will be published separately by the third author. 

First, by \cite[Theorem 5.3]{MR3812195}, one only needs to prove Theorem \ref{thm:transgens} when $G$ is imprimitive with minimal block size $2$, and $n$ has the form $n=2^x3^y5$ with either $y=0$ and $17\le x\le 26$; or $y=1$ and $15\le x\le 35$. Thus, in particular, $G$ may be viewed as a subgroup in a wreath product $2\wr G^{\Sigma}$, where $\Sigma$ is a set of blocks for $G$ of size $2$. It follows that $d(G)\le d_{G^{\Sigma}}(M)+d(G^{\Sigma})$, where $M$ is the intersection of $G$ with the base group of the wreath product, and $d_{G^{\Sigma}}(M)$ is the minimal number of elements required to generate $M$ as a $G^{\Sigma}$-module.
With this reduction in mind, the proof of Theorem \ref{thm:transgens} is comprised of two main ingredients: upper bounds on $d_{G^{\Sigma}}(M)$, and upper bounds on $d(G^{\Sigma})$. We summarise the third author's approach to these two sub-problems in the next few paragraphs.

We note first that the bulk of the proof is taken up with finding upper bounds on $d_{G^{\Sigma}}(M)$. Since $d_{G^{\Sigma}}(M)\le d_H(M)$ for any subgroup $H$ of $G^{\Sigma}$, the strategy of the third author in \cite{MR3812195} in this case involved replacing $G^{\Sigma}$ by a convenient subgroup $H$ of $G^{\Sigma}$, and then deriving upper bounds on $d_H(M)$, usually in terms of the lengths of the $H$-orbits in $\Sigma$. This approach turns out to be particularly fruitful when $H$ is chosen to be a soluble transitive subgroup of $G^{\Sigma}$, whenever such a subgroup exists. When $G^{\Sigma}$ does not contain a soluble transitive subgroup, however, the analysis becomes much more complicated. This led to less sharp bounds, and ultimately, the omitted cases in \cite[Theorem 1.1]{MR3812195}.

The new approach to bounding $d_{G^{\Sigma}}(M)$ involves a careful analysis of the orbit lengths of soluble subgroups in a minimal transitive insoluble subgroup of $G^{\Sigma}$, building on the work in \cite{MR3812195} in the case $n=2^x3$. This analysis, together with upper bounds on $d_H(M)$ (for soluble $H\le G^{\Sigma}$) in terms of the lengths of the $H$-orbits in $\Sigma$, is then used to derive an upper bound for $d_H(M)$. An upper bound for $d_{G^{\Sigma}}(M)$ follows. 

The second sub-problem is to find an upper bound for $d(G^{\Sigma})$.
The group $G^{\Sigma}$ is a transitive permutation group of degree $n/2=2^{x-1}3^y5$, where $n$, $x$, and $y$ are as above. The upper bound $d(G^{\Sigma})\le \frac{c\frac{n}{2}}{\sqrt{\log{\frac{n}{2}}}}$ can be derived by using induction on $n$. However, combining this with the upper bounds on $d(G^{\Sigma})$ detailed in the previous two paragraphs is not enough to prove Theorem \ref{thm:transgens} in all of the required cases. Therefore, a more careful approach is required. This approach was used in the proof of Lemma 5.12 in the third author's paper \cite{MR3812195}. Informally, the idea is as follows. There exists a factorisation $\frac{n}{2}=r_1\hdots r_t$ of $\frac{n}{2}$ such that either
\begin{enumerate}
    \item $d(G^{\Sigma})\le \sum_i=1^{t-1}d(r_i,r_{i+1}\hdots r_t)+\log{r_t}$; or
    \item Either $r_t\le 32$, or $r_t=48$, and $d(G^{\Sigma})\le \sum_i=1^{t-1}d(r_i,r_{i+1}\hdots r_t)+d_{\trans}(r_t)$.
\end{enumerate}
Here, $d_{\trans}(m):=\mathrm{max}\{d(X)\text{ : }X\text{ transitive of degree }m\}$.
If $2\le m\le 32$, or if $m=48$, then we know $d_{\trans}(m)$ precisely, by \cite{MR2455702}, and this paper, respectively.

The function $d(r,s)$ is defined as the maximum of $d_X(K_X(\Delta))$, where
\begin{mylist}
\item[(i)] $X$ runs over the transitive permutation groups of degree $rs$ with minimal block size $r$;
\item[(ii)] $\Delta$ runs over the blocks for $X$ of size $r$;
\item[(iii)] $K_X(\Delta)$ is the kernel of the action of $X$ on $\Delta$; and
\item[(iv)] $d_X(K_X(\Delta))$ is the minimal number of elements required to generated $K_X(\Delta)$ as a normal subgroup of $X$.  
\end{mylist}
Upper bounds on $d(r,s)$ are available from \cite{MR3812195}. Thus, we can find upper bounds on $d(G^{\Sigma})$ by going through all factorisations of $\frac{n}{2}$, and taking the maximum of the bounds coming from (1) and (2) above. These maximums almost always come from (2). Thus, the new result $d_{\trans}(48)=10$ from this paper plays a vital role in deriving upper bounds on $d(G^{\Sigma})$, whence solving the second sub-problem in the proof of Theorem \ref{thm:transgens}.

\subsection*{Acknowledgements}

We would like to thank Michael Giudici for a number of helpful discussions on transitive permutation groups.

\bibliographystyle{plain}

\end{document}